\documentclass{jnmp03c}

\begin{document}

\renewcommand{\evenhead}{J Bystr\"{o}m}
\renewcommand{\oddhead}{Correctors for Some Nonlinear Monotone Operators}

\setcounter{page}{8}
\thispagestyle{empty}

\FistPageHead{1}{\pageref{bystrom-firstpage}--\pageref{bystrom-lastpage}}{Letter}

\copyrightnote{2001}{J Bystr\"{o}m}

\Name{Correctors for Some Nonlinear Monotone Operators}\label{bystrom-firstpage}

\Author{Johan BYSTR\"{O}M}

\Adress{Department of Mathematics, Lule\aa\ University of Technology,
   SE-97187 Lule\aa, Sweden\\
E-mail: johanb@sm.luth.se}

\Date{Received March 3, 2000; Revised June 5, 2000;
Accepted August 28, 2000}

\begin{abstract}
\noindent In this paper we study homogenization of quasi-linear
partial differential equations of the form $-\mbox{div}\left(
a\left( x,x/\varepsilon _h,Du_h\right) \right) =f_h$ on $\Omega $
with Dirichlet boundary conditions. Here the sequence $\left(
\varepsilon _h\right) $ tends to $0$ as $h\rightarrow \infty $ and
the map $a\left( x,y,\xi \right) $ is periodic in~$y,$ monotone in
$\xi $ and satisfies suitable continuity conditions. We prove that
$u_h\rightarrow u$ weakly in $W_0^{1,p}\left( \Omega \right) $ as
$h\rightarrow \infty ,$ where $u$ is the solution of a homogenized
problem of the form $-\mbox{div}\left( b\left( x,Du\right) \right)
=f$ on $\Omega .$ We also derive an explicit expression for the
homogenized operator $b$ and prove some corrector results, i.e. we
find $\left( P_h\right) $ such that $Du_h-P_h\left( Du\right)
\rightarrow 0$ in $L^p\left( \Omega , \mathbf{R}^n\right)$.
\end{abstract}

\section{Introduction}

In mathematical models of microscopically non-homogeneous media, various
local characteristics are usually described by functions of the form
$a\left( x/\varepsilon _h\right) $ where $\varepsilon _h>0$ is a small
parameter. The function $a\left( x\right) $ can be periodic or belong to
some other specific class. To compute the properties of a micro
non-homogeneous medium is an extremely difficult task, since the
coefficients are rapidly oscillating functions. One way to attack the
problem is to apply asymptotic analysis to the problems of microlevel
non-homogeneous media, which immediately leads to the concept of
homogenization. When the parameter $\varepsilon _h$ is very small, the
heterogeneous medium will act as a homogeneous medium. To characterize this
homogeneous medium is one of the main tasks in the homogenization theory.
For more information concerning the homogenization theory, the reader is
referred to \cite{Bens,Jikov} and~\cite{Persson}.

In this paper we consider the homogenization problem for monotone operators
and the local behavior of the solutions. Monotone operators are very
important in the study of nonlinear partial differential equations. The
problem we study here can be used to model different nonlinear stationary
conservation laws, e.g. stationary temperature distribution. For a more
detailed discussion concerning different applications, see~\cite{Zeid1}.

We will study the limit behavior of the sequence of solutions $\left(
u_h\right) $ of the Dirichlet boundary value problem
\[
\left\{
\begin{array}{l}
\ds -\mbox{div}\left( a\left( x,\frac x{\varepsilon
_h},Du_h\right) \right) =f_h\quad \mbox{on}\quad \Omega ,
\vspace{3mm}\\
u_h\in W_0^{1,p}\left( \Omega \right) ,
\end{array}
\right.
\]
where $f_h\rightarrow f$ in $W^{-1,q}\left( \Omega \right) $ as $\varepsilon
_h\rightarrow 0$. Moreover, the map $a\left( x,y,\xi \right) $ is defined on
$\Omega \times \mathbf{R}^n\times \mathbf{R}^n$ and is assumed to be
periodic in $y$, continuous in $\xi $ and monotone in $\xi $. We also need
some continuity restriction in the first variable of $a\left( x,y,\xi
\right) $. We will consider two different cases, namely when $a\left(
x,y,\xi \right) $ is of the form
\[
a\left( x,y,\xi \right) =\sum_{i=1}^N\chi _{\Omega _i}\left( x\right)
a_i\left( y,\xi \right) ,
\]
and when $a\left( x,y,\xi \right) $ satisfies that
\[
\left| a\left( x_1,y,\xi \right) -a\left( x_2,y,\xi \right) \right| ^q\leq
\omega \left( \left| x_1-x_2\right| \right) \left( 1+\left| \xi \right|
\right) ^p,
\]
where $\omega :\mathbf{R}\rightarrow \mathbf{R}$ is continuous, increasing
and $\omega \left( 0\right) =0.$ In both cases we will prove that
$u_h\rightarrow u$ weakly in $W^{1,p}\left( \Omega \right) $ and that $u$ is
the solution of the \textit{homogenized} problem
\[
\left\{
\begin{array}{l}
-\mbox{div}\left( b\left( x,Du\right) \right) =f\quad \mbox{on}\quad \Omega ,
\vspace{2mm}\\
u\in W_0^{1,p}\left( \Omega \right) .
\end{array}
\right.
\]
We will prove that the operator $b$ has the same structure properties as $a$
and is given by
\[
b\left( x,\xi \right) =\int_Ya\left( x,y,\xi +Dv^{\xi ,x}\left( y\right)
\right) \,dy,
\]
where $v^{\xi ,x}$ is the solution of the \textit{cell-problem}
\begin{equation}
\left\{
\begin{array}{l}
-\mbox{div}\left( a\left( x,y,\xi +Dv^{\xi ,x}\left( y\right) \right)
\right) =0\quad \mbox{on}\quad Y,
\vspace{2mm}\\
v^{\xi ,x}\in W_{\square }^{1,p}\left( Y\right) .
\end{array}
\right.  \label{cellproblem0}
\end{equation}
Here $Y$ is a periodic cell and $W_{\square }^{1,p}\left( Y\right) $ is the
subset of $W^{1,p}\left( Y\right) $ such that $u$ has mean value $0$ and $u$
is $Y$-periodic. The corresponding weak formulation of the cell problem is
\begin{equation}
\left\{
\begin{array}{l}
\ds \int_Y\left\langle a\left( x,y,\xi +Dv^{\xi ,x}\left( y\right) \right)
,Dw\right\rangle dy=0\quad \mbox{for every}\quad w\in W_{\square }^{1,p}\left(
Y\right) ,
\vspace{2mm}\\
v^{\xi ,x}\in W_{\square }^{1,p}\left( Y\right) .
\end{array}
\right.  \label{cellvar}
\end{equation}
The homogenization problem described above with $p=2$ was studied in \cite
{Wall1}. Others have investigated the case where we have no dependence in $x$,
that is, when $a$ is of the form $a\left( x,y,\xi \right) =a\left( y,\xi
\right) $. Here we mention~\cite{Fusco} and~\cite{Chiado1} where the
problems corresponding to single valued and multi valued operators were
studied. Moreover, the almost periodic case was treated in~\cite{Brai2}.

The weak convergence of $u_h$ to $u$ in $W^{1,p}\left( \Omega \right) $
implies that $u_h-u\rightarrow 0$ in $L^p\left( \Omega \right) ,$ but in
general, we only have that $Du_h-Du\rightarrow 0$ weakly in $L^p\left(
\Omega ,\mathbf{R}^n\right) $. However, we will prove that it is possible to
express $Du_h$ in terms of $Du$ up to a remainder which converges strongly
in $L^p\left( \Omega ,\mathbf{R}^n\right) $. This is done by constructing a
family of correctors $P_h\left( x,\xi ,t\right) $, defined by
\begin{equation}
P_h\left( x,\xi ,t\right) =P\left( \frac x{\varepsilon _h},\xi ,t\right)
=\xi +Dv^{\xi ,t}\left( \frac x{\varepsilon _h}\right) .  \label{pdef}
\end{equation}

Let $\left( M_h\right) $ be a family of linear operators converging to the
identity map on $L^p\left( \Omega ,\mathbf{R}^n\right) $ such that~$M_hf$ is
a step function for every $f\in L^p\left( \Omega ,\mathbf{R}^n\right) $.
Moreover, let $\gamma _h$ be a step function approximating the identity map
on $\Omega $. We will show that
\[
Du_h-P_h\left( x,M_hDu,\gamma _h\right) \rightarrow 0\quad
\mbox{in}\quad L^p\left(
\Omega ,\mathbf{R}^n\right) .
\]

The problem of finding correctors has been studied by many authors, see e.g.
\cite{Dal1} where single valued monotone operators of the form $-\mbox{div}
\left( a\left( \frac x{\varepsilon _h},Du_h\right) \right) $ were
considered and~\cite{Brai1} where the corresponding almost periodic case was
considered.

\section{Preliminaries and notation}

Let $\Omega $ be a open bounded subset of $\mathbf{R}^n$, $\left| E\right| $
denote the Lebesgue measure in $\mathbf{R}^n$ and $\left\langle \cdot ,\cdot
\right\rangle $ denote the Euclidean scalar product on $\mathbf{R}^n$.
Moreover, if $X$ is a Banach space, we let $X^{*}$ denote its dual space and
$\left\langle \cdot \mid \cdot \right\rangle $ denote the canonical pairing
over $X^{*}\times X$.

Let $\left\{ \Omega _i\subset \Omega :i=1,\ldots ,N\right\} $ be a family of
disjoint open sets such that $\left| \Omega \backslash \cup _{i=1}^N\Omega
_i\right| =0$ and $\left| \partial \Omega _i\right| =0$. Let $\left(
\varepsilon _h\right) $ be a decreasing sequence of real numbers such that
$\varepsilon _h\rightarrow 0$ as $h\rightarrow \infty $. Furthermore,
$Y=\left( 0,1\right) ^n$ is the unit cube in $\mathbf{R}^n$ and we put
$Y_h^j=\varepsilon _h\left( j+Y\right) $, where $j\in \mathbf{Z}^n$, i.e. the
translated image of $\varepsilon _hY$ by the vector $\varepsilon _hj$. We
also define the following index sets:
\[
\ba{l} \ds J_h =\left\{ j\in \mathbf{Z}^n:\overline{Y}_h^j\subset
\Omega \right\}, \qquad J_h^i=\left\{ j\in
\mathbf{Z}^n:\overline{Y}_h^j\subset \Omega _i\right\},
\vspace{3mm}\\
\ds B_h^i =\left\{ j\in \mathbf{Z}^n:\overline{Y}_h^j\cap \Omega _i\neq
\emptyset ,\;\overline{Y}_h^j\backslash \Omega _i\neq \emptyset \right\} .
\ea
\]
Moreover, we define $\Omega _i^h=\cup _{j\in J_h^i}\overline{Y}_h^j$ and
$F_i^h=\cup _{j\in B_h^i}Y_h^j$.

In a corresponding way let $\left\{ \Omega _i^k\subset \Omega :i\in
I_k\right\} $ denote a family of disjoint open sets with diameter less than
$\frac 1k$ such that $\left| \Omega \backslash \cup _{i\in I_k}\Omega
_i^k\right| =0$ and $\left| \partial \Omega _i^k\right| =0$. We also define
the following index sets:
\[
\ba{l}
\ds J_h^{i,k} =\left\{ j\in \mathbf{Z}^n:\overline{Y}_h^j\subset \Omega
_i^k\right\} ,
\vspace{3mm}\\
\ds B_h^{i,k} =\left\{ j\in \mathbf{Z}^n:\overline{Y}_h^j\cap \Omega _i^k\neq
\emptyset ,\;\overline{Y}_h^j\backslash \Omega _i^k\neq \emptyset \right\} .
\ea
\]
Let $\Omega _i^{k,h}=\cup _{j\in J_h^{i,k}}\overline{Y}_h^j$ and
$F_i^{k,h}=\cup _{j\in B_h^{i,k}}Y_h^j$.

Corresponding to $f\in L^p\left( \Omega ,\mathbf{R}^n\right) $ we define the
function $M_hf:\mathbf{R}^n\rightarrow \mathbf{R}^n$ by
\[
\left( M_hf\right) \left( x\right) =\sum_{j\in J_h}\chi _{Y_h^j}\left(
x\right) \xi _h^j,
\]
where $\xi _h^j=\frac 1{\left| Y_h^j\right| }\int_{Y_h^j}f\,dx$ and $\chi _E$
is the characteristic function of the set $E$ (in order to define $\xi _h^j$
for all $j\in \mathbf{Z}^n,$ we treat $f$ as $f=0$ outside $\Omega $). It is
well known that
\begin{equation}
M_hf\rightarrow f\quad \mbox{in}\quad L^p(\Omega ,\mathbf{R}^n),  \label{approx}
\end{equation}
see \cite[p.~129]{Royden}. We also define the step function $\gamma
_h:\Omega \rightarrow \Omega $ by
\begin{equation}
\gamma _h\left( x\right) =\sum_{j\in J_h}\chi _{Y_h^j}\left( x\right) x_h^j,
\label{gdef}
\end{equation}
where $x_h^j\in Y_h^j.$ Finally, $C$ will denote a positive constant that
may differ from one place to an other.

Let $a:\Omega \times \mathbf{R}^n\times \mathbf{R}^n\rightarrow
\mathbf{R}^n$ be a function such that $a\left( x,\cdot ,\xi
\right) $ is Lebesgue measurable and $Y$-periodic for $x\in \Omega
$ and $\xi \in \mathbf{R}^n.$ Let $p$ be a real constant
$1<p<\infty $ and let $q$ be its dual exponent, $\frac 1p+\frac
1q=1.$ We also assume that $a$ satisfies the following continuity
and monotonicity conditions: There exists two constants $c_1,$
$c_2>0,$ and two constants $\alpha $ and~$\beta ,$ with $0\leq
\alpha \leq \min \left( 1,p-1\right) $ and $\max \left( p,2\right)
\leq \beta <\infty $ such that \be \left| a\left( x,y,\xi
_1\right) -a\left( x,y,\xi _2\right) \right| \leq c_1\left(
1+\left| \xi _1\right| +\left| \xi _2\right| \right) ^{p-1-\alpha
}\left| \xi _1-\xi _2\right| ^\alpha , \label{acont} \ee \be
\left\langle a\left( x,y,\xi _1\right) -a\left( x,y,\xi _2\right)
,\xi _1-\xi _2\right\rangle \geq c_2\left( 1+\left| \xi _1\right|
+\left| \xi _2\right| \right) ^{p-\beta }\left| \xi _1-\xi
_2\right| ^\beta , \label{amon} \ee for $x\in \Omega,$ a.e. $y\in
\mathbf{R}^n$ and every $\xi \in \mathbf{R}^n$. Moreover, we
assume that
\begin{equation}
a\left( x,y,0\right) =0,  \label{axy0}
\end{equation}
for $x\in \Omega,$ a.e. $y\in \mathbf{R}^n$. Let $\left( f_h\right) $ be a
sequence in $W^{-1,q}\left( \Omega \right) $ that converges to $f$.

\medskip

\noindent
{\bf Remark 1.}
We will use these continuity and monotonicity conditions to show theorems
and properties. However, we concentrate on showing the non-trivial cases,
for instance when $\beta \neq p,$ and omit the simple ones, in this case
when $\beta =p.$

\medskip

The solution $v^{\xi ,x}$ of the cell-problem (\ref{cellproblem0})
can be extended by periodicity to an element in $W_{\rm
loc}^{1,p}\left( \mathbf{R}^n\right)$, still denoted by $v^{\xi
,x},$ and
\begin{equation}
\int_{\mathbf{R}^n}\left\langle a\left( x,y,\xi +Dv^{\xi ,x}\left( y\right)
\right) ,D\phi \left( y\right) \right\rangle \,dy=0 \quad \mbox{for every}\quad \phi
\in C_0^\infty \left( \mathbf{R}^n\right) .  \label{divwhk}
\end{equation}

\section{Some useful lemmas}

The following lemma, see e.g. \cite{Murat}, is fundamental to the
homogenization theory.

\medskip

\noindent {\bf Lemma 1 (Compensated compactness).} {\it Let
$1<p<\infty $. Moreover, let $\left( v_h\right) $ be a sequence in
$L^q\left( \Omega , \mathbf{R}^n\right) $ which converges weakly
to $v$, $\left( -\mbox{\rm div}\, v_h\right) $ converges to
$-\mbox{\rm div}\, v$ in $W^{-1,q}(\Omega )$ and let $\left(
u_h\right) $ be a sequence which converges weakly to $u$ in
$W^{1,p}\left( \Omega \right) .$ Then
\[
\int_\Omega \left\langle v_h,Du_h\right\rangle \phi \,dx\rightarrow
\int_\Omega \left\langle v,Du\right\rangle \phi \,dx,
\]
for every $\phi \in C_0^\infty \left( \Omega \right) .$}

\medskip

We will also use the following estimates, which are proved in~\cite{Bys}.

\medskip

\noindent
{\bf Lemma 2.}
{\it Let $a$ satisfy (\ref{acont}), (\ref{amon}) and (\ref{axy0}). Then the
following inequalities hold:
\be
(a) \quad
\left| a\left( x,y,\xi \right) \right| \leq c_a\left( 1+\left| \xi \right|
^{p-1}\right) ,  \label{ineqa}
\ee
\be
(b)\quad
\left| \xi \right| ^p\leq c_b\left( 1+\left\langle a\left( x,y,\xi \right)
,\xi \right\rangle \right) ,  \label{ineqb}
\ee
\be
(c) \quad
\int_Y\left| \xi +Dv^{\xi ,x_i}\right| ^pdy\leq c_c\left( 1+\left| \xi
\right| ^p\right) .  \label{ineqc}
\ee}

\noindent
{\bf Lemma 3.} {\it For every $\xi _1,\xi _2\in \mathbf{R}^n$ we have
\begin{equation}
\int_Y\left| \xi _1+Dv^{\xi _1,x}-\xi _2-Dv^{\xi _2,x}\right| ^pdy\leq
C\left( 1+\left| \xi _1\right| ^p+\left| \xi _2\right| ^p\right) ^{\frac{
\beta -\alpha -1}{\beta -\alpha }}\left| \xi _1-\xi _2\right| ^{\frac
p{\beta -\alpha }}.
\label{dalmasoineq}
\end{equation}}

\section{Some homogenization results}

Let $a\left( x,y,\xi \right) $ satisfy one of the conditions
\begin{enumerate}
\topsep0mm
\partopsep0mm
\parsep0mm
\itemsep0mm
\item $a$ is of the form
\begin{equation}
a\left( x,y,\xi \right) =\sum_{i=1}^N\chi _{\Omega _i}\left( x\right)
a_i\left( y,\xi \right) .  \label{eq14}
\end{equation}
\item There exists a function $\omega :\mathbf{R\rightarrow R}$ that is
continuous, increasing and $\omega \left( 0\right) =0$ such that
\begin{equation}
\left| a\left( x_1,y,\xi \right) -a\left( x_2,y,\xi \right) \right| ^q\leq
\omega \left( \left| x_1-x_2\right| \right) \left( 1+\left| \xi \right|
\right) ^p  \label{eq15}
\end{equation}
for $x_1,$ $x_2\in \Omega ,$ a.e. $y\in \mathbf{R}$ and every $\xi \in
\mathbf{R}^n.$
\end{enumerate}

Now we consider the weak Dirichlet boundary value problems (one for each
choice of~$h$):
\begin{equation}
\left\{
\begin{array}{l}
\ds \int_\Omega \left\langle a\left( x,\frac x{\varepsilon _h},Du_h\right)
,D\phi \right\rangle dx=\left\langle f_h\mid \phi \right\rangle \quad \mbox{for
every}\quad \phi \in W_0^{1,p}\left( \Omega \right) ,
\vspace{3mm}\\
\ds u_h\in W_0^{1,p}\left( \Omega \right) .
\end{array}
\right.  \label{eq16}
\end{equation}
By a standard result in the existence theory for boundary value problems
defined by monotone operators, these problems have a unique solution for
each~$h,$ see e.g.~\cite{Zeid2}.

We let $\phi =u_h$ in (\ref{eq16}) and use H\"{o}lder's reversed inequality,
(\ref{amon}), (\ref{axy0}) and Poincare's inequality. This implies that
\be\label{eq17}
\ba{l}
\ds
C\left( \int_\Omega \left( 1+\left| Du_h\right| \right) ^pdx\right)^{\frac{p-\beta }p}
\left( \int_\Omega \left| Du_h\right| ^pdx\right)
^{\frac \beta p}\leq c_2\int_\Omega \left( 1+\left| Du_h\right| \right)
^{p-\beta }\left| Du_h\right|^\beta dx
\vspace{3mm}\\
\ds \qquad \qquad
\leq \int_\Omega \left\langle a\left( x,\frac x{\varepsilon
_h},Du_h\right) ,Du_h\right\rangle dx=\left\langle f_h\mid u_h\right\rangle
\vspace{3mm}\\
\ds \qquad \qquad \leq \left\| f_h\right\| _{W^{-1,q}}\left\| u_h\right\| _{W_0^{1,p}}\leq
C\left\| Du_h\right\| _{L^p\left( \Omega ,\mathbf{R}^n\right) },
\ea\hspace{-10mm}
\ee
where $C$ does not depend on $h.$ Now if $\left\| Du_h\right\| _{L^p\left(
\Omega ,\mathbf{R}^n\right) }^p<\left| \Omega \right| ,$ then clearly
$\left\| u_h\right\| _{W_0^{1,p}\left( \Omega \right) }\leq C$ by Poincare's
inequality. Hence assume that $\left\| Du_h\right\| _{L^p\left( \Omega ,
\mathbf{R}^n\right) }^p\geq \left| \Omega \right| .$ But then we have
by~(\ref{eq17}) that
\[
\ba{l}
\ds 2^{p-\beta }\int_\Omega \left| Du_h\right|^pdx=\left( 2^{p-1}\int_\Omega
2\left| Du_h\right| ^pdx\right) ^{\frac{p-\beta }p}\left( \int_\Omega
\left| Du_h\right| ^pdx\right) ^{\frac \beta p}
\vspace{3mm}\\
\ds \qquad
\leq \left( \int_\Omega 2^{p-1}\left( 1+\left| Du_h\right| ^p\right)
dx\right) ^{\frac{p-\beta }p}\left( \int_\Omega \left| Du_h\right|^pdx\right)^{\frac \beta p}
\leq C\left( \int_\Omega \left| Du_h\right|^pdx\right)^{\frac 1p},
\ea
\]
that is, $\left\| Du_h\right\| _{L^p\left( \Omega ,\mathbf{R}^n\right) }\leq
C.$ According to Poincare's inequality we thus have that
\linebreak $\left\|
u_h\right\| _{W_0^{1,p}\left( \Omega \right) }\leq C.$ Summing up, we have
that
\begin{equation}
\left\| u_h\right\| _{W_0^{1,p}\left( \Omega \right) }\leq C.
\label{ubounded}
\end{equation}
Since $u_h$ is bounded in $W_0^{1,p}\left( \Omega \right) $, there exists a
subsequence $\left( h^{\prime }\right) $ such that
\begin{equation}
u_{h^{\prime }}\rightarrow u_{*}\quad \mbox{weakly in} \quad
W_0^{1,p}\left( \Omega \right) .  \label{uconv}
\end{equation}
The following theorems will show that $u_{*}$ satisfy an equation of the
same type as those which are satisfied by $u_h.$

\medskip

\noindent
{\bf Theorem 4.} {\it Let $a$ satisfy (\ref{acont}), (\ref{amon}), (\ref{axy0}) and
(\ref{eq14}). Let $\left( u_h\right) $ be solutions of (\ref{eq16}). Then we
have that
\be
\ba{l}
u_h\rightarrow u\quad \mbox{weakly in} \quad W_0^{1,p}\left( \Omega \right) ,
\vspace{3mm}\\
\ds a\left( x,\frac x{\varepsilon _h},Du_h\right) \rightarrow b\left(
x,Du\right) \quad \mbox{weakly in} \quad  L^q\left( \Omega ,\mathbf{R}^n\right) ,
\ea \label{ahbounded}
\ee
where $u$ is the unique solution of the homogenized problem
\begin{equation}
\left\{
\begin{array}{l}
\ds \int_\Omega \left\langle b\left( x,Du\right) ,D\phi
\right\rangle dx=\left\langle f\mid \phi \right\rangle \quad
\mbox{for every}\quad \phi \in W_0^{1,p}\left( \Omega \right) ,
\vspace{3mm}\\
u\in W_0^{1,p}\left( \Omega \right) .
\ea
\right.  \label{eq19}
\end{equation}
The operator $b:\Omega \times \mathbf{R}^n\rightarrow \mathbf{R}^n$ is
defined a.e. as
\[
b\left( x,\xi \right) =\sum_{i=1}^N\chi _{\Omega _i}\left( x\right)
\int_Ya_i\left( y,\xi +Dv^{\xi ,x_i}\left( y\right) \right)
dy=\sum_{i=1}^N\chi _{\Omega _i}\left( x\right) b_i\left( \xi \right) ,
\]
where $x_i\in \Omega _i,$ $b_i\left( \xi \right) =\int_Ya_i\left( y,\xi
+Dv^{\xi ,x_i}\left( y\right) \right) dy$ and $v^{\xi ,x_i}$ is the unique
solution of the cell problem
\begin{equation}
\left\{
\begin{array}{l}
\ds \int_Y\left\langle a_i\left( y,\xi +Dv^{\xi ,x_i}\left(
y\right) \right) ,D\phi \left( y\right) \right\rangle dy=0\quad
\mbox{for every} \quad \phi \in W_{\Box }^{1,p}\left( Y\right) ,
\vspace{3mm}\\
\ds v^{\xi ,x_i}\in W_{\Box }^{1,p}\left( Y\right) .
\end{array}
\right.  \label{eq21}
\end{equation}}

\noindent
{\bf Remark 2.}
An equivalent formulation of the equations (\ref{eq19}) and (\ref{eq21})
above can be given in the following unified manner
\[
\left\{
\begin{array}{l}
\ds \int_\Omega \int_Y\left\langle a\left( x,y,Du\left( x\right)
+Dv\left( x,y\right) \right) ,D\overline{u}\left( x\right)
+D\overline{v}\left( x,y\right) \right\rangle dxdy=\int_\Omega
f\overline{u}\, dx,
\vspace{3mm}\\
\ds \left( u,v\right) \in W_0^{1,p}\left( \Omega \right) \times L^q\left( \Omega
;W_{\Box }^{1,p}\left( Y\right) \right) ,
\end{array}
\right.
\]
for all $\left( \overline{u},\overline{v}\right) \in W_0^{1,p}\left( \Omega
\right) \times L^q\left( \Omega ;W_{\Box }^{1,p}\left( Y\right) \right) .$
This formulation often occurs in the notion of two-scale limit introduced in
e.g.~\cite{All}.

\medskip

\noindent
{\bf Proof.} We have shown that $u_{h^{\prime }}\rightarrow u_{*}$ weakly in
$W_0^{1,p}\left( \Omega \right) $ for a subsequence $\left( h^{\prime
}\right) $ since $u_h$ is bounded in $W_0^{1,p}\left( \Omega \right) .$ We
define
\[
\psi _{h^{\prime }}^i=a_i\left( \frac x{\varepsilon _{h^{\prime}}},Du_{h^{\prime }}\right) .
\]
Then according to (\ref{acont}), (\ref{axy0}), H\"{o}lder's inequality,
Poincare's inequality and (\ref{ubounded}), we find that
\[
\ba{l}
\ds \int_{\Omega _i}\left| \psi _{h^{\prime }}^i\right|^qdx=\int_{\Omega
_i}\left| a_i\left( \frac x{\varepsilon _{h^{\prime }}},Du_{h^{\prime
}}\right) \right|^qdx\leq c_1^q\int_{\Omega _i}\left( 1+\left|
Du_{h^{\prime }}\right| \right)^{q\left( p-1-\alpha \right) }\left|
Du_{h^{\prime }}\right| ^{\alpha q}dx
\vspace{3mm}\\
\ds \phantom{\int_{\Omega _i}\left| \psi _{h^{\prime }}^i\right|^qdx=}
 \leq C\left( \int_{\Omega _i}\left( 1+\left| Du_{h^{\prime }}\right|
\right) ^pdx\right) ^{\frac{p-1-\alpha }{p-1}}\left( \int_{\Omega
_i}\left| Du_{h^{\prime }}\right| ^pdx\right) ^{\frac \alpha {p-1}}
\vspace{3mm}\\
\ds \phantom{\int_{\Omega _i}\left| \psi _{h^{\prime }}^i\right|^qdx=}
\leq C\int_{\Omega _i}\left( 1+\left| Du_{h^{\prime }}\right| \right)^pdx\leq
C,
\ea
\]
that is, $\psi _{h^{\prime }}^i$ is bounded in $L^q\left( \Omega _i,\mathbf{R}^n\right) .$
Hence there is a subsequence $\left( h^{\prime \prime }\right)
$ of $\left( h^{\prime }\right) $ such that
\[
\psi _{h^{\prime \prime }}^i\rightarrow \psi _{*}^i\quad \mbox{weakly in}
\quad L^q\left( \Omega _i,\mathbf{R}^n\right) .
\]
From our original problem (\ref{eq16}) we conclude that
\[
\left\{
\begin{array}{l}
\ds \sum\limits_{i=1}^N\int_{\Omega _i}\left\langle a_i\left( \frac
x{\varepsilon _{h^{\prime \prime }}},Du_{h^{\prime \prime }}\right) ,D\phi
\right\rangle dx=\left\langle f_{h^{\prime \prime }}\mid \phi \right\rangle
\quad \mbox{for every}\quad \phi \in W_0^{1,p}\left( \Omega \right) ,
\vspace{3mm}\\
\ds u_{h^{\prime \prime }}\in W_0^{1,p}\left( \Omega \right) .
\end{array}
\right.
\]
In the limit we have
\[
\sum_{i=1}^N\int_{\Omega _i}\left\langle \psi _{*}^i,D\phi \right\rangle
dx=\left\langle f\mid \phi \right\rangle \quad \mbox{for every}\quad \phi \in
W_0^{1,p}\left( \Omega \right) .
\]
If we could show that
\begin{equation}
\psi _{*}^i=b_i\left( Du_{*}\right) \quad \mbox{for a.e.}\quad x\in \Omega _i,
\label{eq24b}
\end{equation}
then it follows by uniqueness of the homogenized problem (\ref{eq19}) that
$u_{*}=u$. Let
\[
w_h^{\xi ,x_i}\left( x\right) =\left\langle \xi ,x\right\rangle +\varepsilon
_hv^{\xi ,x_i}\left( \frac x{\varepsilon _h}\right)
\]
for a.e. $x\in \mathbf{R}^n.$ By the monotonicity of $a_i$ we have for a fix
$\xi $ that
\[
\int_{\Omega _i}\!\left\langle \! a_i\left( \frac x{\varepsilon _{h^{\prime
\prime }}},Du_{h^{\prime \prime }}\left( x\right) \right) -a_i\left( \frac
x{\varepsilon _{h^{\prime \prime }}},Dw_{h^{\prime \prime }}^{\xi
,x_i}\left( x\right) \right) ,Du_{h^{\prime \prime }}\left( x\right)
-Dw_{h^{\prime \prime }}^{\xi ,x_i}\left( x\right) \!\right\rangle \phi \left(
x\right) dx\geq 0
\]
for every $\phi \in C_0^\infty \left( \Omega _i\right) $, $\phi \geq 0.$ We
now note that by (\ref{eq16}) it also holds that
\[
\int_{\Omega _i}\left\langle a_i\left( \frac x{\varepsilon _{h^{\prime
\prime }}},Du_{h^{\prime \prime }}\left( x\right) \right) ,D\phi
\right\rangle dx=\left\langle f_{h^{\prime \prime }}\mid \phi \right\rangle
\quad \mbox{for every}\quad \phi \in W_0^{1,p}\left( \Omega _i\right) .
\]
Since $f_h\rightarrow f$ in $W^{-1,q}\left( \Omega \right) ,$ this implies
that
\begin{equation}
-\mbox{div}\left( a_i\left( \frac x{\varepsilon _{h^{\prime \prime
}}},Du_{h^{\prime \prime }}\left( x\right) \right) \right) =-\mbox{div}\,
\psi _{h^{\prime \prime }}^i\rightarrow -\mbox{div}\,\psi _{*}^i
\quad \mbox{in} \quad W^{-1,q}\left( \Omega _i\right) .  \label{eq27}
\end{equation}
We also have that
\begin{equation}
-\mbox{div}\left( a_i\left( \frac x{\varepsilon _{h^{\prime \prime
}}},Dw_{h^{\prime \prime }}^{\xi ,x_i}\left( x\right) \right) \right) =0
\quad \mbox{on} \quad \Omega _i.  \label{eq28}
\end{equation}
By the compensated compactness lemma (Lemma~1) we
then get in the limit
\[
\int_{\Omega _i}\left\langle \psi _{*}^i-b_i\left( \xi \right) ,Du_{*}\left(
x\right) -\xi \right\rangle \phi \left( x\right) dx\geq 0
\]
for every $\phi \in C_0^\infty \left( \Omega _i\right)$, $\phi \geq 0.$
Hence for our fixed $\xi \in \mathbf{R}^n$ we have that
\begin{equation}
\left\langle \psi _{*}^i-b_i\left( \xi \right) ,Du_{*}\left( x\right) -\xi
\right\rangle \geq 0\quad \mbox{for a.e.} \quad x\in \Omega _i. \label{eq30}
\end{equation}
In particular, if $\left( \xi _m\right) $ is a countable dense set in
$\mathbf{R}^n$, then (\ref{eq30}) implies that
\begin{equation}
\left\langle \psi _{*}^i-b_i\left( \xi _m\right) ,Du_{*}\left( x\right) -\xi
_m\right\rangle \geq 0 \quad \mbox{for a.e.}\quad x\in \Omega _i.  \label{eq31}
\end{equation}
By the continuity of $b$ (see Remark~3) it follows that
\begin{equation}
\left\langle \psi _{*}^i-b_i\left( \xi \right) ,Du_{*}\left(
x\right) -\xi \right\rangle \geq 0\quad \mbox{for a.e.} \quad x\in
\Omega _i \quad \mbox{and for every} \quad \xi \in \mathbf{R}^n.
\label{eq32}
\end{equation}
Since $b$ is monotone and continuous (see Remark~3) we have that
$b$ is maximal monotone and hence (\ref{eq24b}) follows. We have proved the
theorem up to a subsequence $\left( u_{h^{\prime \prime }}\right) $ of
$\left( u_h\right) .$ By uniqueness of the solution to the homogenized
equation it follows that this holds for the whole sequence.
\hfill \rule{3mm}{3mm}

\medskip

\noindent
{\bf Theorem 5.} {\it Let $a$ satisfy (\ref{acont}), (\ref{amon}), (\ref{axy0}) and
(\ref{eq15}). Let $\left( u_h\right) $ be solutions of (\ref{eq16}). Then we
have that
\[
\ba{l} u_h\rightarrow u \quad \mbox{weakly in}\quad
W_0^{1,p}\left(\Omega \right),
\vspace{3mm}\\
\ds a\left( x,\frac x{\varepsilon _h},Du_h\right) \rightarrow
b\left( x,Du\right) \quad \mbox{weakly in}\quad L^q\left( \Omega
,\mathbf{R}^n\right) , \ea
\]
where $u$ is the unique solution of the homogenized problem
\[
\left\{
\begin{array}{l}
\ds \int_\Omega \left\langle b\left( x,Du\right) ,D\phi \right\rangle
dx=\left\langle f\mid \phi \right\rangle \quad \mbox{for every} \quad \phi \in
W_0^{1,p}\left( \Omega \right) ,
\vspace{3mm}\\
\ds u\in W_0^{1,p}\left( \Omega \right) .
\end{array}
\right.
\]
The operator $b:\Omega \times \mathbf{R}^n\rightarrow \mathbf{R}^n$ is
defined a.e. as
\[
b\left( x,\xi \right) =\int_Ya\left( x,y,\xi +Dv^{\xi ,x}\left( y\right)
\right) dy,
\]
where $v^{\xi ,x}$ is the unique solution of the cell problem
\begin{equation}
\left\{
\begin{array}{l}
\ds \int_Y\left\langle a\left( x,y,\xi +Dv^{\xi ,x}\left( y\right) \right)
,D\phi \left( y\right) \right\rangle dy=0 \quad \mbox{for every} \quad \phi \in W_{\Box
}^{1,p}\left( Y\right) ,
\vspace{3mm}\\
\ds v^{\xi ,x}\in W_{\Box }^{1,p}\left( Y\right) .
\end{array}
\right.  \label{eq35}
\end{equation}}

Before we prove this theorem, we make some definitions that will be useful
in the proof. Define the function
\[
a^k\left( x,y,\xi \right) =\sum_{i\in I_k}\chi _{\Omega _i^k}\left( x\right)
a\left( x_i^k,y,\xi \right) ,
\]
where $x_i^k\in \Omega _i^k.$ Consider the boundary value problems
\begin{equation}
\left\{
\begin{array}{l}
\ds \int_\Omega \left\langle a^k\left( x,\frac x{\varepsilon _h},Du_h^k\right)
,D\phi \right\rangle dx=\left\langle f_h\mid \phi \right\rangle
\quad  \mbox{for every}\quad \phi \in W_0^{1,p}\left( \Omega \right) ,
\vspace{3mm}\\
u_h^k\in W_0^{1,p}\left( \Omega \right) .
\end{array}
\right.  \label{eq37}
\end{equation}
The conditions for Theorem~4 are satisfied and the theorem implies
that there exists a $u_{*}^k$ such that
\[
u_h^k\rightarrow u_{*}^k\quad \mbox{weakly in} \quad W_0^{1,p}\left( \Omega \right)
\quad \mbox{as} \quad h\rightarrow \infty ,
\]
where $u_{*}^k$ is the unique solution of
\begin{equation}
\left\{
\begin{array}{l}
\ds \int_\Omega \left\langle b^k\left( x,Du_{*}^k\right) ,D\phi \right\rangle
dx=\left\langle f\mid \phi \right\rangle \quad \mbox{for every} \quad \phi \in
W_0^{1,p}\left( \Omega \right) ,
\vspace{3mm}\\
\ds u_{*}^k\in W_0^{1,p}\left( \Omega \right) .
\end{array}
\right.  \label{eq38}
\end{equation}
Here
\[
b^k\left( x,\xi \right) =\sum_{i\in I_k}\chi _{\Omega _i^k}\left( x\right)
\int_Ya\left( x_i^k,y,\xi +Dv^{\xi ,x_i^k}\left( y\right) \right)
dy=\sum_{i\in I_k}\chi _{\Omega _i^k}\left( x\right) b\left( x_i^k,\xi
\right) ,
\]
and $v^{\xi ,x_i^k}$ is the solution of
\[
\left\{
\begin{array}{l}
\ds \int_Y\left\langle a\left( x_i^k,y,\xi +Dv^{\xi ,x_i^k}\left( y\right)
\right) ,D\phi \left( y\right) \right\rangle dy=0\quad \mbox{for every} \quad \phi \in
W_{\Box }^{1,p}\left( Y\right) ,
\vspace{3mm}\\
v^{\xi ,x_i^k}\in W_{\Box }^{1,p}\left( Y\right) .
\end{array}
\right.
\]

\noindent
{\bf Proof.} First we prove that $u_h\rightarrow u$ weakly in $W_0^{1,p}\left( \Omega
\right) .$ If $g\in W^{-1,q}\left( \Omega \right) ,$ we have that
\[
\ba{l}
\ds \lim\limits_{h\rightarrow \infty }\left\langle g\mid u_h-u\right\rangle
=\lim\limits_{k\rightarrow \infty }\lim\limits_{h\rightarrow \infty
}\left\langle g\mid u_h-u\right\rangle
\vspace{3mm}\\
\ds \qquad \qquad
=\lim\limits_{k\rightarrow \infty
}\lim\limits_{h\rightarrow \infty }\left\langle g\mid
u_h-u_h^k+u_h^k-u_{*}^k+u_{*}^k-u\right\rangle
\vspace{3mm}\\
\ds \qquad \qquad
\leq \lim\limits_{k\rightarrow \infty }\lim\limits_{h\rightarrow \infty
}\left\| g\right\| _{W^{-1,q}\left( \Omega \right) }\left\|
u_h-u_h^k\right\| _{W_0^{1,p}\left( \Omega \right)
}+\lim\limits_{k\rightarrow \infty }\lim\limits_{h\rightarrow \infty
}\left\langle g\mid u_h^k-u_{*}^k\right\rangle
\vspace{3mm}\\
\ds \qquad \qquad
+\lim\limits_{k\rightarrow \infty }\left\| g\right\| _{W^{-1,q}\left(
\Omega \right) }\left\| u_{*}^k-u\right\| _{W_0^{1,p}\left( \Omega \right)}.
\ea
\]
We need to show that all terms on the right hand side are zero.

\medskip

\noindent
\textbf{Term 1.}  Let us prove that
\begin{equation}
\lim_{k\rightarrow \infty }\lim_{h\rightarrow \infty }\left\|
u_h-u_h^k\right\| _{W_0^{1,p}\left( \Omega \right) }=0.  \label{eq42}
\end{equation}
By the definition we have that
\[
\ba{l}
\ds \int_\Omega \left\langle a^k\left( x,\frac x{\varepsilon
_h},Du_h^k\right) ,D\phi \right\rangle dx=\left\langle f_h\mid \phi
\right\rangle \quad \mbox{for every}\quad \phi \in W_0^{1,p}\left( \Omega \right) ,
\vspace{3mm}\\
\ds \int_\Omega \left\langle a\left( x,\frac x{\varepsilon _h},Du_h\right)
,D\phi \right\rangle dx=\left\langle f_h\mid \phi \right\rangle \quad \mbox{for every}
\quad \phi \in W_0^{1,p}\left( \Omega \right) .
\ea
\]
This implies that we with $\phi =u_h^k-u_h$ have
\[
\ba{l} \ds \int_\Omega \left\langle a^k\left( x,\frac
x{\varepsilon _h},Du_h^k\right) -a^k\left( x,\frac x{\varepsilon
_h},Du_h\right) ,Du_h^k-Du_h\right\rangle dx
\vspace{3mm}\\
\ds \qquad =\int_\Omega \left\langle a\left( x,\frac x{\varepsilon _h},Du_h\right)
-a^k\left( x,\frac x{\varepsilon _h},Du_h\right) ,Du_h^k-Du_h\right\rangle
dx.
\ea
\]
By the monotonicity of $a$ (\ref{amon}), Schwarz' inequality and
H\"{o}lder's inequality it follows that
\[
\ba{l}
\ds c_2\int_\Omega \left( 1+\left| Du_h^k\right| +\left| Du_h\right| \right)
^{p-\beta }\left| Du_h^k-Du_h\right| ^\beta dx
\vspace{3mm}\\
\ds \qquad \leq \int_\Omega \left\langle a^k\left( x,\frac x{\varepsilon
_h},Du_h^k\right) -a^k\left( x,\frac x{\varepsilon _h},Du_h\right)
,Du_h^k-Du_h\right\rangle dx
\vspace{3mm}\\
\ds \qquad =\int_\Omega \left\langle a\left( x,\frac x{\varepsilon _h},Du_h\right)
-a^k\left( x,\frac x{\varepsilon _h},Du_h\right) ,Du_h^k-Du_h\right\rangle
dx
\vspace{3mm}\\
\ds \qquad \leq \left( \int_\Omega \left| a\left( x,\frac
x{\varepsilon _h},Du_h\right) -a^k\left( x,\frac x{\varepsilon
_h},Du_h\right) \right| ^qdx\right) ^{\frac 1q}\left( \int_\Omega
\left| Du_h^k-Du_h\right| ^pdx\right) ^{\frac 1p}. \ea
\]
Now since $\left\| Du_h^k\right\| _{L^p\left( \Omega ,\mathbf{R}^n\right)
}^p\leq C$ and $\left\| Du_h\right\| _{L^p\left( \Omega ,\mathbf{R}^n\right)
}^p\leq C,$ we have, according to H\"{o}lder's reversed inequality, that
\be
\ba{l}
\ds C\left( \int_\Omega \left| Du_h^k-Du_h\right| ^pdx\right) ^{\frac \beta
p}
\vspace{3mm}\\
\ds \qquad \leq c_2\left( \int_\Omega \left( 1+\left| Du_h^k\right| +\left|
Du_h\right| \right) ^pdx\right) ^{\frac{p-\beta }p}\left( \int_\Omega
\left| Du_h^k-Du_h\right| ^pdx\right) ^{\frac \beta p}
\vspace{3mm}\\
\ds \qquad
\leq c_2\int_\Omega \left( 1+\left| Du_h^k\right| +\left| Du_h\right|
\right) ^{p-\beta }\left| Du_h^k-Du_h\right| ^\beta dx.
\ea
\label{eq47}
\ee

Hence we see that
\[
\left\| Du_h^k-Du_h\right\| _{L^p\left( \Omega ,\mathbf{R}^n\right) }^p\leq
C\left( \int_\Omega \left| a\left( x,\frac x{\varepsilon _h},Du_h\right)
-a^k\left( x,\frac x{\varepsilon _h},Du_h\right) \right| ^qdx\right)^{\frac{p-1}{\beta -1}}.
\]
Thus in view of the continuity condition (\ref{eq15}) it yields that
\begin{equation}
\ba{l}
\ds \left\| Du_h^k-Du_h\right\| _{L^p\left( \Omega ,\mathbf{R}^n\right) }^p
\vspace{3mm}\\
\ds \qquad \leq C\left( \omega \left( \frac 1k\right) \int_\Omega \left( 1+\left|
Du_h\right| \right) ^pdx\right) ^{\frac{p-1}{\beta -1}}\leq C\left( \omega
\left( \frac 1k\right) \right) ^{\frac{p-1}{\beta -1}},
\ea\label{eq43}
\end{equation}
where we in the last inequality have used the fact that there exists a
constant $C$ independent of $h$ such that $\left\| Du_h\right\| _{L^p\left(
\Omega ,\mathbf{R}^n\right) }\leq C.$ Since $\left\| D\cdot \right\|
_{L^p\left( \Omega ,\mathbf{R}^n\right) }$ is an equivalent norm on
$W_0^{1,p}\left( \Omega \right) $, we have that
\[
\left\| u_h-u_h^k\right\| _{W_0^{1,p}\left( \Omega \right) }\rightarrow 0
\]
as $k\rightarrow \infty $ uniformly in $h.$ This means that we can change
order in the limit process in~(\ref{eq42}) and (\ref{eq42}) follows by
taking (\ref{eq43}) into account.

\medskip

\noindent
\textbf{Term 2.}  We observe that
\[
\lim_{k\rightarrow \infty }\lim_{h\rightarrow \infty }\left\langle g\mid
u_h^k-u_{*}^k\right\rangle =0
\]
as a direct consequence of Theorem~4.

\medskip

\noindent
\textbf{Term 3.}  Let us prove that
\begin{equation}
\lim_{k\rightarrow \infty }\left\| u_{*}^k-u\right\| _{W_0^{1,p}\left(
\Omega \right) }=0.  \label{eq44}
\end{equation}
By the definition we have that
\[
\ba{l} \ds \int_\Omega \left\langle b^k\left( x,Du_{*}^k\right)
,D\phi \right\rangle dx=\left\langle f\mid \phi \right\rangle
\quad \mbox{for every}\quad \phi \in W_0^{1,p}\left( \Omega
\right) ,
\vspace{3mm}\\
\ds \int_\Omega \left\langle b\left( x,Du\right) ,D\phi \right\rangle
dx=\left\langle f\mid \phi \right\rangle \quad \mbox{for every} \quad \phi \in
W_0^{1,p}\left( \Omega \right) .
\ea
\]
Hence
\[
\int_\Omega \left\langle b^k\left( x,Du_{*}^k\right) -b^k\left( x,Du\right)
,D\phi \right\rangle dx=\int_\Omega \left\langle b\left( x,Du\right)
-b^k\left( x,Du\right) ,D\phi \right\rangle dx
\]
for every $\phi \in W_0^{1,p}\left( \Omega \right) .$ Now let $\phi
=u_{*}^k-u$ and use the strict monotonicity of $b^k$ (see Remark~3)
on the left hand side and use Schwarz' and H\"{o}lder's
inequalities on the right hand side. We get
\[
\ba{l}
\ds c_2\int_\Omega \left( 1+\left| Du_{*}^k\right| +\left| Du\right| \right)
^{p-\beta }\left| Du_{*}^k-Du\right| ^\beta dx
\vspace{3mm}\\
\ds \qquad \leq \int_\Omega \left\langle b^k\left( x,Du_{*}^k\right) -b^k\left(
x,Du\right) ,Du_{*}^k-Du\right\rangle dx
\vspace{3mm}\\
\ds \qquad =\int_\Omega \left\langle b\left( x,Du\right) -b^k\left( x,Du\right)
,Du_{*}^k-Du\right\rangle dx
\vspace{3mm}\\
\ds \qquad \leq \left( \int_\Omega \left| b\left( x,Du\right) -b^k\left( x,Du\right)
\right| ^qdx\right) ^{\frac 1q}\left( \int_\Omega \left|
Du_{*}^k-Du\right| ^pdx\right) ^{\frac 1p}.
\ea
\]

Inequality (\ref{eq47}) with $Du_{*}^k$ and $Du$ instead of $Du_h^k$ and
$Du_h,$ respectively, then yields that
\begin{equation}
\left\| Du_{*}^k-Du\right\| _{L^p\left( \Omega ,\mathbf{R}^n\right) }^p\leq
C\left( \omega \left( \frac 1k\right) \int_\Omega \left( 1+\left|
Du\right| \right) ^pdx\right) ^{\frac{p-1}{\beta -1}}\leq C\left( \omega
\left( \frac 1k\right) \right) ^{\frac{p-1}{\beta -1}}.  \label{eq45}\hspace{-5mm}
\end{equation}

The right hand side tends to $0$ as $k\rightarrow \infty .$ We now obtain
(\ref{eq44}) by noting that $\left\| D\cdot \right\| _{L^p\left( \Omega ,
\mathbf{R}^n\right) }$ is an equivalent norm on $W_0^{1,p}\left( \Omega
\right) $.

\medskip

Next we prove that $a\left( x,\frac x{\varepsilon _h},Du_h\right)
\rightarrow b\left( x,Du\right) $ weakly in $L^q\left( \Omega
,\mathbf{R} ^n\right) .$ Now if $g\in \left( L^q\left( \Omega
,\mathbf{R}^n\right) \right)^{*}$, then
\[
\ba{l}
\ds \lim\limits_{h\rightarrow \infty }\! \left\langle g\mid a\left( x,\frac
x{\varepsilon _h},Du_h\right) -b\left( x,Du\right) \right\rangle
=\lim\limits_{k\rightarrow \infty }\lim\limits_{h\rightarrow \infty
}\! \left\langle g\mid a\left( x,\frac x{\varepsilon _h},Du_h\right) -b\left(
x,Du\right) \right\rangle
\vspace{3mm}\\
\ds \qquad =\lim\limits_{k\rightarrow \infty }\lim\limits_{h\rightarrow \infty
}\left\langle g\mid a\left( x,\frac x{\varepsilon _h},Du_h\right)
-a^k\left( x,\frac x{\varepsilon _h},Du_h^k\right) \right\rangle
\vspace{3mm}\\
\ds \qquad +\lim\limits_{k\rightarrow \infty }\lim\limits_{h\rightarrow \infty
}\left\langle g\mid a^k\left( x,\frac x{\varepsilon _h},Du_h^k\right)
-b^k\left( x,Du_{*}^k\right) \right\rangle
\vspace{3mm}\\
\ds \qquad +\lim\limits_{k\rightarrow \infty
}\lim\limits_{h\rightarrow \infty }\left\langle g\mid b^k\left(
x,Du_{*}^k\right) -b\left( x,Du\right) \right\rangle
\vspace{3mm}\\
\ds \qquad \leq \lim\limits_{k\rightarrow \infty }\lim\limits_{h\rightarrow \infty
}\left\| g\right\| _{\left( L^q\left( \Omega ,\mathbf{R}^n\right) \right)
^{*}}\left\| a\left( x,\frac x{\varepsilon _h},Du_h\right) -a^k\left(
x,\frac x{\varepsilon _h},Du_h^k\right) \right\| _{L^q\left( \Omega ,
\mathbf{R}^n\right) }
\vspace{3mm}\\
\ds \qquad +\lim\limits_{k\rightarrow \infty }\lim\limits_{h\rightarrow \infty
}\left\langle g\mid a^k\left( x,\frac x{\varepsilon _h},Du_h^k\right)
-b^k\left( x,Du_{*}^k\right) \right\rangle
\vspace{3mm}\\
\ds \qquad +\lim\limits_{k\rightarrow \infty }\left\| g\right\| _{\left( L^q\left(
\Omega ,\mathbf{R}^n\right) \right) ^{*}}\left\| b^k\left( x,Du_{*}^k\right)
-b\left( x,Du\right) \right\| _{L^q\left( \Omega ,\mathbf{R}^n\right) }.
\ea\!\!
\]
We now prove that the three terms on the right hand side are zero.

\strut\hfill

\noindent
\textbf{Term 1.}  Let us show that
\begin{equation}
\lim_{k\rightarrow \infty }\lim_{h\rightarrow \infty }\left\| a\left(
x,\frac x{\varepsilon _h},Du_h\right) -a^k\left( x,\frac x{\varepsilon
_h},Du_h^k\right) \right\| _{L^q\left( \Omega ,\mathbf{R}^n\right) }=0.
\label{eq46}
\end{equation}
By Minkowski's inequality we have that
\[
\ba{l}
\ds \int_\Omega \left| a^k\left( x,\frac x{\varepsilon _h},Du_h^k\right)
-a\left( x,\frac x{\varepsilon _h},Du_h\right) \right| ^qdx
\vspace{3mm}\\
\ds \qquad
\leq C\int_\Omega \left| a^k\left( x,\frac x{\varepsilon
_h},Du_h^k\right) -a^k\left( x,\frac x{\varepsilon _h},Du_h\right) \right|
^qdx
\vspace{3mm}\\
\ds \qquad +C\int_\Omega \left| a^k\left( x,\frac x{\varepsilon _h},Du_h\right)
-a\left( x,\frac x{\varepsilon _h},Du_h\right) \right| ^qdx.
\ea
\]
The second integral on the right hand side is bounded by the continuity
condition (\ref{eq15}) since
\begin{equation}
\ba{l}
\ds \int_\Omega \left| a^k\left( x,\frac x{\varepsilon _h},Du_h\right)
-a\left( x,\frac x{\varepsilon _h},Du_h\right) \right| ^qdx
\vspace{3mm}\\
\ds \qquad \leq \omega \left( \frac 1k\right) \int_\Omega \left( 1+\left|
Du_h\right| \right) ^pdx\leq C\omega \left( \frac 1k\right) ,
\ea\label{eq49}
\end{equation}
where the last inequality follows since $\left( Du_h\right) $ is bounded in
$L^p\left( \Omega ,\mathbf{R}^n\right) .$ For the first integral on the right
hand side we use the continuity restriction~(\ref{acont}), H\"{o}lder's
inequality and~(\ref{eq43}) to derive that
\be \label{eq50}
\ba{l}
\ds \int_\Omega \left| a^k\left( x,\frac x{\varepsilon _h},Du_h^k\right)
-a^k\left( x,\frac x{\varepsilon _h},Du_h\right) \right| ^qdx
\vspace{3mm}\\
\ds \qquad \leq c_1^q\int_\Omega \left( 1+\left| Du_h^k\right| +\left| Du_h\right|
\right) ^{q(p-1-\alpha )}\left| Du_h^k-Du_h\right| ^{\alpha q}dx
\vspace{3mm}\\
\ds \qquad \leq C\left( \int_\Omega \left( 1+\left| Du_h^k\right| +\left| Du_h\right|
\right) ^pdx\right)^{\frac{p-1-\alpha }{p-1}}\left( \int_\Omega \left|
Du_h^k-Du_h\right| ^pdx\right) ^{\frac \alpha {p-1}}
\vspace{3mm}\\
\ds \qquad \leq C\left( \int_\Omega \left| Du_h^k-Du_h\right| ^pdx\right) ^{\frac
\alpha {p-1}}\leq C\left( \omega \left( \frac 1k\right) \right) ^{\frac
\alpha {\beta -1}}
\ea
\ee
since $\left( Du_h\right) $ and $\left( Du_h^k\right) $ are bounded in
$L^p\left( \Omega ,\mathbf{R}^n\right) .$ Hence by (\ref{eq49}) and
(\ref{eq50}) we have that
\begin{equation}
\int_\Omega \left| a^k\left( x,\frac x{\varepsilon _h},Du_h^k\right)
-a\left( x,\frac x{\varepsilon _h},Du_h\right) \right| ^qdx\leq C\left(
\omega \left( \frac 1k\right) \right) ^{\frac \alpha {\beta -1}}+C\omega
\left( \frac 1k\right) .  \label{eq54}
\end{equation}
By the properties of $\omega ,$ it follows that
\[
\left\| a\left( x,\frac x{\varepsilon _h},Du_h\right) -a^k\left( x,\frac
x{\varepsilon _h},Du_h^k\right) \right\| _{L^q\left( \Omega ,\mathbf{R}^n\right) }
\rightarrow 0
\]
as $k\rightarrow \infty $ uniformly in $h$. This implies that we may change
the order in the limit process in (\ref{eq46}) and (\ref{eq46}) follows.

\medskip

\noindent
\textbf{Term 2.}  We immediately see that
\[
\lim_{k\rightarrow \infty }\lim_{h\rightarrow \infty }\left\langle g\mid
a^k\left( x,\frac x{\varepsilon _h},Du_h^k\right) -a^k\left( x,\frac
x{\varepsilon _h},Du_{*}^k\right) \right\rangle =0
\]
as a direct consequence of Theorem~4.

\medskip

\noindent
\textbf{Term 3.}  Let us show that
\begin{equation}
\lim_{k\rightarrow \infty }\left\| b^k\left( x,Du_{*}^k\right) -b\left(
x,Du\right) \right\| _{L^q\left( \Omega ,\mathbf{R}^n\right) }=0.
\label{eq55}
\end{equation}
We note that
\[
\ba{l}
\ds \int_\Omega \left| b^k\left( x,Du_{*}^k\right) -b\left( x,Du\right)
\right| ^qdx
\vspace{3mm}\\
\ds \qquad =\int_\Omega \left| b^k\left( x,Du_{*}^k\right) -b^k\left( x,Du\right)
+b^k\left( x,Du\right) -b\left( x,Du\right) \right| ^qdx
\vspace{3mm}\\
\ds \qquad \leq C\left( \int_\Omega \left| b^k\left( x,Du_{*}^k\right) -b^k\left(
x,Du\right) \right| ^qdx+\int_\Omega \left| b^k\left( x,Du\right) -b\left(
x,Du\right) \right| ^qdx\right) .
\ea
\]
By applying the continuity conditions in Remark~3, Theorem~6
 and H\"{o}lder's inequality, we see that
\[
\ba{l}
\ds \int_\Omega \left| b^k\left( x,Du_{*}^k\right) -b\left( x,Du\right)
\right| ^qdx
\vspace{3mm}\\
\ds \qquad \leq C\left( \int_\Omega \left( 1+\left| Du_{*}^k\right| +\left| Du\right|
\right) ^pdx\right) ^{\frac{p-1-\gamma }{p-1}}\left( \int_\Omega \left|
Du_{*}^k-Du\right| ^pdx\right) ^{\frac \gamma {p-1}}
\vspace{3mm}\\
\ds \qquad +C\widetilde{\omega }\left( \frac 1k\right) \int_\Omega 1+\left|
Du\right| ^pdx\leq C\left( \int_\Omega \left| Du_{*}^k-Du\right|
^pdx\right) ^{\frac \gamma {p-1}}+C\widetilde{\omega }\left( \frac
1k\right) .
\ea
\]
Now (\ref{eq55}) follows by taking (\ref{eq45}) into account.\hfill\rule{3mm}{3mm}

\section{Properties of the homogenized operator}

In this section we prove some properties of the homogenized operator. In
particular, these properties imply the existence and uniqueness of the
solution of the homogenized problem.

\medskip

\noindent
{\bf Theorem 6.} {\it Let $b$ be the homogenized operator defined in Theorem~5.
Then
\begin{enumerate}
\topsep0mm
\partopsep0mm
\parsep0mm
\itemsep0mm
\item[(a)]  $b\left( \cdot ,\xi \right) $ satisfies the continuity condition
\[
\left| b\left( x_1,\xi \right) -b\left( x_2,\xi \right) \right| ^q\leq
\widetilde{\omega }\left( \left| x_1-x_2\right| \right) \left( 1+\left| \xi
\right| \right) ^p,
\]
where $\widetilde{\omega }:\mathbf{R\rightarrow R}$ is a function that is
continuous, increasing and $\widetilde{\omega }\left( 0\right) =0$.
\item[(b)]  $b\left( x,\cdot \right) $ is strictly monotone. In particular,
we have that
\[
\left\langle b\left( x,\xi _1\right) -b\left( x,\xi _2\right) ,\xi _1-\xi
_2\right\rangle \geq \widetilde{c}_2\left( 1+\left| \xi _1\right| +\left|
\xi _2\right| \right) ^{p-\beta }\left| \xi _1-\xi _2\right| ^\beta
\]
for every $\xi _1,\xi _2\in \mathbf{R}^n.$
\item[(c)]  $b\left( x,\cdot \right) $ is continuous. In particular, we have
for $\gamma =\frac \alpha {\beta -\alpha }$ that
\[
\left| b\left( x,\xi _1\right) -b\left( x,\xi _2\right) \right| \leq
\widetilde{c}_1\left( 1+\left| \xi _1\right| +\left| \xi _2\right| \right)
^{p-1-\gamma }\left| \xi _1-\xi _2\right| ^\gamma
\]
for every $\xi _1,\xi _2\in \mathbf{R}^n.$
\item[(d)]  $b\left( x,0\right) =0$ for $x\in \Omega .$
\end{enumerate}}

\noindent
{\bf Proof.}

$(a)$  By the definition of $b$ we have that
\[
\ba{l}
\ds \left| b\left( x_1,\xi \right) -b\left( x_2,\xi \right) \right| ^q=\left|
\int_Ya\left( x_1,y,\xi +Dv^{\xi ,x_1}\right) dy-\int_Ya\left( x_2,y,\xi
+Dv^{\xi ,x_2}\right) dx\right| ^q
\vspace{3mm}\\
\ds
\qquad \leq C\int_Y\left| a\left( x_1,y,\xi +Dv^{\xi ,x_1}\right) -a\left(
x_2,y,\xi +Dv^{\xi ,x_1}\right) \right| ^qdy
\vspace{3mm}\\
\ds \qquad +C\int_Y\left| a\left( x_2,y,\xi +Dv^{\xi ,x_1}\right) -a\left( x_2,y,\xi
+Dv^{\xi ,x_2}\right) \right| ^qdy,
\ea
\]
where the last inequality follows from Jensen's inequality. By the
continuity conditions~(\ref{eq15}) and (\ref{acont}), it then
follows that \be \label{eq51} \ba{l} \ds \left| b\left( x_1,\xi
\right) -b\left( x_2,\xi \right) \right| ^q\leq C\omega \left(
\left| x_1-x_2\right| \right) \int_Y\left( 1+\left| \xi +Dv^{\xi
,x_1}\right| \right) ^pdy
\vspace{3mm}\\
\ds \qquad
+C\left( \int_Y\left( 1+\left| \xi +Dv^{\xi ,x_1}\right| ^p+\left| \xi
+Dv^{\xi ,x_2}\right| ^p\right) dy\right) ^{\frac{p-1-\alpha }{p-1}}
\vspace{3mm}\\
\ds \qquad \times \left( \int_Y\left| Dv^{\xi ,x_1}-Dv^{\xi ,x_2}\right| ^pdy\right)
^{\frac \alpha {p-1}}.
\ea
\ee
We now study the two terms above separately. For the first term we have by
(\ref{ineqc}) that
\begin{equation}
\int_Y\left( 1+\left| \xi +Dv^{\xi ,x_i}\right| \right) ^pdy\leq
C\left( 1+\left| \xi \right| ^p\right). \label{eq52}
\end{equation}
For the second term we have by the definition that
\[
\ba{l}
\ds \int_Y\left\langle a\left( x_1,y,\xi +Dv^{\xi ,x_1}\left( y\right) \right)
,D\phi \right\rangle dy=0 \quad \mbox{for every} \quad
\phi \in W_{\square }^{1,p}\left(\Omega \right) ,
\vspace{3mm}\\
\ds \int_Y\left\langle a\left( x_2,y,\xi +Dv^{\xi ,x_2}\left( y\right) \right)
,D\phi \right\rangle dy=0\quad \mbox{for every} \quad \phi \in W_{\square }^{1,p}
\left(\Omega \right) .
\ea
\]
This implies that
\[
\ba{l}
\ds \int_Y\left\langle a\left( x_1,y,\xi +Dv^{\xi ,x_1}\left( y\right) \right)
-a\left( x_1,y,\xi +Dv^{\xi ,x_2}\left( y\right) \right) ,D\phi
\right\rangle dy
\vspace{3mm}\\
\ds \qquad =\int_Y\left\langle a\left( x_2,y,\xi +Dv^{\xi ,x_2}\left( y\right)
\right) -a\left( x_1,y,\xi +Dv^{\xi ,x_2}\left( y\right) \right) ,D\phi
\right\rangle dy
\ea
\]
for every $\phi \in W_{\square }^{1,p}\left( \Omega \right) .$ In
particular, for $\phi =v^{\xi ,x_1}-v^{\xi ,x_2}$ we have that
\[
\ba{l}
\ds \int_Y\left\langle a\left( x_1,y,\xi +Dv^{\xi ,x_1}\left( y\right) \right)
-a\left( x_1,y,\xi +Dv^{\xi ,x_2}\left( y\right) \right) ,Dv^{\xi
,x_1}-Dv^{\xi ,x_2}\right\rangle dy
\vspace{3mm}\\
\ds \qquad =\int_Y\left\langle a\left( x_2,y,\xi +Dv^{\xi ,x_2}\left( y\right)
\right) -a\left( x_1,y,\xi +Dv^{\xi ,x_2}\left( y\right) \right) ,Dv^{\xi
,x_1}-Dv^{\xi ,x_2}\right\rangle dy.
\ea
\]
We apply (\ref{amon}) on the left hand side and Schwarz' and H\"{o}lder's
inequalities together with~(\ref{eq15}) and~(\ref{eq52}) on the right hand
side. This yields
\[
\ba{l}
\ds c_2\int_Y\left( 1+\left| \xi +Dv^{\xi ,x_1}\right| +\left| \xi +Dv^{\xi
,x_2}\right| \right) ^{p-\beta }\left| Dv^{\xi ,x_1}-Dv^{\xi ,x_2}\right|
^\beta dy
\vspace{3mm}\\
\ds \qquad \leq \int_Y\left\langle a\left( x_1,y,\xi +Dv^{\xi ,x_1}\left( y\right)
\right) -a\left( x_1,y,\xi +Dv^{\xi ,x_2}\left( y\right) \right) ,Dv^{\xi
,x_1}-Dv^{\xi ,x_2}\right\rangle dy
\vspace{3mm}\\
\ds \qquad =\int_Y\left\langle a\left( x_2,y,\xi +Dv^{\xi ,x_2}\left( y\right)
\right) -a\left( x_1,y,\xi +Dv^{\xi ,x_2}\left( y\right) \right) ,Dv^{\xi
,x_1}-Dv^{\xi ,x_2}\right\rangle dy
\vspace{3mm}\\
\ds \qquad \leq \left( \int_Y\left| a\left( x_2,y,\xi +Dv^{\xi ,x_2}\left( y\right)
\right) -a\left( x_1,y,\xi +Dv^{\xi ,x_2}\left( y\right) \right) \right|
^qdy\right) ^{\frac 1q}
\vspace{3mm}\\
\ds \qquad
\times \left( \int_Y\left| Dv^{\xi ,x_1}-Dv^{\xi ,x_2}\right| ^pdy\right)
^{\frac 1p}
\vspace{3mm}\\
\ds \qquad \leq \omega \left( \left| x_1-x_2\right| \right) ^{\frac 1q}C\left(
1+\left| \xi \right| ^p\right) ^{\frac 1q}\left( \int_Y\left| Dv^{\xi
,x_1}-Dv^{\xi ,x_2}\right| ^pdy\right) ^{\frac 1p}.
\ea
\]
The reversed H\"{o}lder inequality then ensures that
\[
\ba{l}
\ds \left( \int_Y\left( 1+\left| \xi +Dv^{\xi ,x_1}\right| +\left| \xi
+Dv^{\xi ,x_2}\right| \right) ^pdy\right) ^{\frac{p-\beta }p}\left(
\int_Y\left| Dv^{\xi ,x_1}-Dv^{\xi ,x_2}\right| ^pdy\right) ^{\frac \beta p}
\vspace{3mm}\\
\ds \qquad \leq \int_Y\left( 1+\left| \xi +Dv^{\xi ,x_1}\right| +\left| \xi +Dv^{\xi
,x_2}\right| \right) ^{p-\beta }\left| Dv^{\xi ,x_1}-Dv^{\xi ,x_2}\right|
^\beta dy,
\ea
\]
that is,
\be\label{eq53}
\ba{l}
\ds \left( \int_Y\left| Dv^{\xi ,x_1}-Dv^{\xi ,x_2}\right| ^pdy\right)
^{\frac \alpha {p-1}}\leq C\omega \left( \left| x_1-x_2\right| \right)
^{\frac \alpha {\beta -1}}\left( 1+\left| \xi \right| ^p\right) ^{\frac
\alpha {\beta -1}}
\vspace{3mm}\\
\ds \qquad \times \left( \int_Y\left( 1+\left| \xi +Dv^{\xi ,x_1}\right| +\left| \xi
+Dv^{\xi ,x_2}\right| \right) ^pdy\right) ^{\frac{\alpha \left( \beta
-p\right) }{\left( \beta -1\right) \left( p-1\right) }}
\vspace{3mm}\\
\ds \qquad \leq C\omega \left( \left| x_1-x_2\right| \right) ^{\frac \alpha {\beta
-1}}\left( 1+\left| \xi \right| ^p\right) ^{\frac \alpha {\beta -1}}\left(
1+\left| \xi \right| ^p\right) ^{\frac{\alpha \left( \beta -p\right) }
{\left( \beta -1\right) \left( p-1\right) }}
\vspace{3mm}\\
\ds \qquad =C\omega \left( \left| x_1-x_2\right| \right) ^{\frac \alpha {\beta
-1}}\left( 1+\left| \xi \right| ^p\right) ^{\frac \alpha {p-1}},
\ea
\ee
where the last inequality follows from (\ref{eq52}). By combining (\ref{eq52})
and~(\ref{eq53}), we can write~(\ref{eq51}) as
\be
\ba{l}
\ds \left| b\left( x_1,\xi \right) -b\left( x_2,\xi \right) \right| ^q\leq
C\omega \left( \left| x_1-x_2\right| \right) \left( 1+\left| \xi \right|
^p\right)
\vspace{3mm}\\
\ds \qquad+C\left( \int_Y\left( 1+\left| \xi +Dv^{\xi ,x_1}\right| ^p+\left| \xi
+Dv^{\xi ,x_2}\right| ^p\right) dy\right) ^{\frac{p-1-\alpha }{p-1}}
\vspace{3mm}\\
\ds \qquad \times \omega
\left( \left| x_1-x_2\right| \right) ^{\frac \alpha {\beta -1}}\left(
1+\left| \xi \right| ^p\right) ^{\frac \alpha {p-1}}
\vspace{3mm}\\
\ds \qquad \leq C\left( 1+\left| \xi \right| ^p\right) \left( \omega \left( \left|
x_1-x_2\right| \right) +\omega \left( \left| x_1-x_2\right| \right) ^{\frac
\alpha {\beta -1}}\right) =\widetilde{\omega }\left( \left| x_1-x_2\right|
\right) \left( 1+\left| \xi \right| ^p\right) ,
\ea\label{eq56}\hspace{-10mm}
\ee
where $\widetilde{\omega }:\mathbf{R\rightarrow R}$ is a function that is
continuous, increasing and $\widetilde{\omega }\left( 0\right) =0$.

\medskip

$(b)$  Let $\xi _j\in \mathbf{R}^n,$ $j=1,2$ and define for a.e. $x\in \mathbf{R}^n$
\[
w_h^{\xi _j,x}=\left( \xi _j,y\right) +\varepsilon _hv^{\xi _j,x}\left(
\frac y{\varepsilon _h}\right) .
\]
By the periodicity of $v^{\xi _j,x}$ we have that
\be
 w_h^{\xi _j,x}\rightarrow \left( \xi _j,y\right) \quad \mbox{weakly in} \quad
W^{1,p}\left( Y\right) ,
\label{eq57}
\ee
\be
Dw_h^{\xi _j,x}\rightarrow \xi _j \quad \mbox{weakly in} \quad L^p\left( Y,\mathbf{R}
^n\right) ,
\label{eq58}
\ee
\be
a\left( x,\frac y{\varepsilon _h},Dw_h^{\xi _j,x}\right) \rightarrow
b\left( x,\xi _j\right) \quad \mbox{weakly in} \quad L^p\left( Y,\mathbf{R}^n\right) .
\label{eq59}
\ee
Moreover, the monotonicity condition (\ref{amon}) on $a$ implies that
\[
\ba{l}
\ds \int_Y\left\langle a\left( x,\frac y{\varepsilon _h},Dw_h^{\xi
_1,x}\right) -a\left( x,\frac y{\varepsilon _h},Dw_h^{\xi _2,x}\right)
,Dw_h^{\xi _1,x}-Dw_h^{\xi _2,x}\right\rangle \phi \left( y\right) dy
\vspace{3mm}\\
\ds \qquad \geq c_2\int_Y\left( 1+\left| Dw_h^{\xi _1,x}\right| +\left| Dw_h^{\xi
_2,x}\right| \right) ^{p-\beta }\left| Dw_h^{\xi _1,x}-Dw_h^{\xi
_2,x}\right| ^\beta \phi \left( y\right) dy
\ea
\]
for every $\phi \in C_0^\infty \left( \Omega \right) $ such that $\phi \geq
0.$ The reversed H\"{o}lder inequality then yields
\[
\ba{l}
\ds \int_Y\left\langle a\left( x,\frac y{\varepsilon _h},Dw_h^{\xi
_1,x}\right) -a\left( x,\frac y{\varepsilon _h},Dw_h^{\xi _2,x}\right)
,Dw_h^{\xi _1,x}-Dw_h^{\xi _2,x}\right\rangle \phi \left( y\right) dy
\vspace{3mm}\\
\ds \qquad \geq c_2\left( \int_Y\left( 1+\left| Dw_h^{\xi _1,x}\right| +\left|
Dw_h^{\xi _2,x}\right| \right) ^p\phi \left( y\right) dy\right) ^{\frac{p-\beta }p}
\vspace{3mm}\\
\ds \qquad \times \left( \int_Y\left| Dw_h^{\xi _1,x}-Dw_h^{\xi _2,x}\right| ^p\phi
\left( y\right) dy\right) ^{\frac \beta p}.
\ea
\]
We apply $\liminf\limits_{h\rightarrow \infty }$ on both sides of this inequality
and obtain
\[
\ba{l}
\ds \left\langle b\left( x,\xi _1\right) -b\left( x,\xi _2\right) ,\xi _1-\xi
_2\right\rangle \int_Y\phi \left( y\right) dy
=\int_Y\left\langle b\left( x,\xi _1\right) -b\left( x,\xi _2\right) ,\xi
_1-\xi _2\right\rangle \phi \left( y\right) dy
\vspace{3mm}\\
\ds \qquad \geq c_2\liminf\limits_{h\rightarrow \infty }\left( \int_Y\left( 1+\left|
Dw_h^{\xi _1,x}\right| +\left| Dw_h^{\xi _2,x}\right| \right) ^p\phi \left(
y\right) dy\right) ^{\frac{p-\beta }p}
\vspace{3mm}\\
\ds \qquad \times \left( \int_Y\left| Dw_h^{\xi _1,x}-Dw_h^{\xi _2,x}\right| ^p\phi
\left( y\right) dy\right) ^{\frac \beta p}
\ea
\]
\[
\ba{l}
\ds \qquad \geq C\liminf\limits_{h\rightarrow \infty }\left( \int_Y\left( 1+\left|
Dw_h^{\xi _1,x}\right| ^p+\left| Dw_h^{\xi _2,x}\right| ^p\right) \phi
\left( y\right) dy\right) ^{\frac{p-\beta }p}
\vspace{3mm}\\
\ds \qquad \times \liminf\limits_{h\rightarrow \infty }\left( \int_Y\left| Dw_h^{\xi
_1,x}-Dw_h^{\xi _2,x}\right| ^p\phi \left( y\right) dy\right) ^{\frac \beta p}
\vspace{3mm}\\
\ds \qquad \geq C\left( \int_Y\left( 1+\left| \xi _1\right| ^p+\left| \xi _2\right|
^p\right) \phi \left( y\right) dy\right) ^{\frac{p-\beta }p}\left(
\int_Y\left| \xi _1-\xi _2\right| ^p\phi \left( y\right) dy\right) ^{\frac \beta p}
\vspace{3mm}\\
\ds \qquad \geq C\left( 1+\left| \xi _1\right| +\left| \xi _2\right| \right) ^{p-\beta
}\left| \xi _1-\xi _2\right| ^\beta \int_Y\phi \left( y\right) dy,
\ea
\]
for every $\phi \in C_0^\infty \left( \Omega \right) .$ Here we have used
(\ref{eq58}), (\ref{eq59}) and Lemma~1 on the left
hand side, and on the right hand side we have used~(\ref{eq58}) and the fact
that for a weakly convergent sequence $\left( x_n\right) $ converging to $x,$
we have that $\left\| x\right\| \leq \liminf\limits_{h\rightarrow \infty }\left\|
x_n\right\| .$ This implies that
\[
\left\langle b\left( x,\xi _1\right) -b\left( x,\xi _2\right) ,\xi _1-\xi
_2\right\rangle \geq \widetilde{c}_2\left( 1+\left| \xi _1\right| +\left|
\xi _2\right| \right) ^{p-\beta }\left| \xi _1-\xi _2\right| ^\beta .
\]
Since $\xi _1$, $\xi _2$ were arbitrarily chosen, $(b)$ is proved.

\medskip

$(c)$  Fix $\xi _1,\xi _2\in \mathbf{R}^n.$ According to Jensen's
inequality, (\ref{acont}), H\"{o}lder's inequality, (\ref{ineqc})
and~(\ref{dalmasoineq}) we have that
\[
\ba{l}
\ds \left| b\left( x,\xi _1\right) -b\left( x,\xi _2\right) \right| ^q=\left|
\int_Ya\left( x,y,\xi _1+Dv^{\xi _1,x}\right) dy-\int_Ya\left( x,y,\xi
_2+Dv^{\xi _2,x}\right) dy\right| ^q
\vspace{3mm}\\
\ds \qquad \leq \int_Y\left| a\left( x,y,\xi _1+Dv^{\xi _1,x}\right) -a\left( x,y,\xi
_2+Dv^{\xi _2,x}\right) \right| ^qdy
\vspace{3mm}\\
\ds \qquad \leq \int_Yc_1^q\left( 1+\left| \xi _1+Dv^{\xi _1,x}\right| +\left| \xi
_2+Dv^{\xi _2,x}\right| \right) ^{\left( p-1-\alpha \right) q}
\vspace{3mm}\\
\ds \qquad \times \left| \xi
_1+Dv^{\xi _1,x}-\xi _2-Dv^{\xi _2,x}\right| ^{\alpha q}dy
\vspace{3mm}\\
\ds \qquad \leq C\left( \int_Y\left( 1+\left| \xi _1+Dv^{\xi _1,x}\right| +\left| \xi
_2+Dv^{\xi _2,x}\right| \right) ^pdy\right) ^{\frac{p-1-\alpha }{p-1}}
\vspace{3mm}\\
\ds \qquad \times \left( \int_Y\left| \xi _1+Dv^{\xi _1,x}-\xi _2-Dv^{\xi
_2,x}\right| ^pdy\right) ^{\frac \alpha {p-1}}
\vspace{3mm}\\
\ds \qquad \leq C\left( 1+\left| \xi _1\right| ^p+\left| \xi
_2\right| ^p \right)^{\frac{p-1-\alpha }{p-1}}\left( \left(
1+\left| \xi _1\right| ^p+\left| \xi _2\right| ^p\right)
^{\frac{\beta -\alpha -1}{\beta -\alpha }}\left| \xi _1-\xi
_2\right| ^{\frac p{\beta -\alpha }}\right) ^{\frac \alpha {p-1}}
\vspace{3mm}\\
\ds \qquad
\leq C\left( 1+\left| \xi _1\right| +\left| \xi _2\right| \right) ^{q\left(
p-1-\frac \alpha {\beta -a}\right) }\left| \xi _1-\xi _2\right| ^{q\frac
\alpha {\beta -a}}.
\ea
\]
Hence
\[
\left| b\left( x,\xi _1\right) -b\left( x,\xi _2\right) \right| \leq
\widetilde{c}_1\left( 1+\left| \xi _1\right| +\left| \xi _2\right| \right)
^{p-1-\gamma }\left| \xi _1-\xi _2\right| ^\gamma ,
\]
where
\[
\gamma =\frac \alpha {\beta -a}.
\]

(d)  Since $a\left( x,y,0\right) =0$ we have that the solution of the
cell problem (\ref{eq35}) corresponding to $\xi =0$ is $v^{0,x}=0.$ This
implies that
\[
b\left( x,0\right) =\int_Ya\left( x,y,0\right) dy=0.
\]
The proof is complete.\hfill \rule{3mm}{3mm}

\medskip

\noindent
{\bf Remark 3.} By using similar arguments as in the proof above we find that
(b), (c) and~(d) holds, up to boundaries, also for the homogenized operator
$b$ in Theorem~4.

\section{Some corrector results}

We have proved in both Theorem~4 and Theorem~5 that we for
the corresponding solutions have that $u_h-u$ converges to $0$ weakly in
$W_0^{1,p}\left( \Omega \right) $. By the Rellich imbedding theorem, we have
that $u_h-u$ converges to $0$ in $L^p\left( \Omega \right) $. In general, we
do not have strong convergence of $Du_h-Du$ to $0$ in $L^p\left( \Omega ,
\mathbf{R}^n\right) $. However, we will prove that it is possible to express
$Du_h$ in terms of $Du$, up to a remainder which converges to $0$ in
$L^p\left( \Omega ,\mathbf{R}^n\right) $.

\medskip

\noindent
{\bf Theorem 7.} {\it Let $u$ and $u_h$ be defined as in Theorem~4 and
let $P_h$ be given by (\ref{pdef}). Then
\[
Du_h-\sum_{i=1}^N\chi _{\Omega _i}\left( x\right) P_h\left(
x,M_hDu,x_i\right) \rightarrow 0\quad \mbox{in} \quad L^p\left( \Omega ,
\mathbf{R}^n\right) .
\]}

\noindent
{\bf Proof.} In \cite{Dal1} the case $N=1$ was considered and in~\cite{Wall1}
 the case $p=2$ was considered.
By using these ideas and making the necessary
adjustments the proof follows. For the details, see~\cite{Bys}.
\hfill \rule{3mm}{3mm}

\medskip

\noindent
{\bf Theorem 8.} {\it Let $u$ and $u_h$ be defined as in Theorem~5.
Moreover, let $P_h$ be given by (\ref{pdef}) and $\gamma _h$ by (\ref{gdef}). Then
\[
Du_h-P_h\left( x,M_hDu,\gamma _h\right) \rightarrow 0 \quad
\mbox{in} \quad L^p\left( \Omega ,\mathbf{R}^n\right) .
\]}

\noindent
{\bf Proof.} We have that
\be\label{equ1thconcor}
\ba{l}
\ds \left\| Du_h-M_hDu-Dv^{M_hDu,\gamma _h}\left( \frac x{\varepsilon
_h}\right) \right\| _{L^p\left( \Omega ,\mathbf{R}^n\right) }
\leq \left\| Du_h-Du_h^k\right\| _{L^p\left( \Omega ,\mathbf{R}^n\right)}
\vspace{3mm}\\
\ds
+\left\| Du_h^k-M_hDu_{*}^k-\sum\limits_{i\in I_k}\chi _{\Omega
_i^k}\left( x\right) Dv^{M_hDu_{*}^k,x_i^k}\left( \frac x{\varepsilon
_h}\right) \right\| _{L^p\left( \Omega ,\mathbf{R}^n\right) }
\vspace{3mm}\\
\ds  +\left\| M_hDu_{*}^k+\sum\limits_{i\in I_k}\chi _{\Omega _i^k}\left(
x\right) Dv^{M_hDu_{*}^k,x_i^k}\left( \frac x{\varepsilon _h}\right)
-M_hDu-Dv^{M_hDu,\gamma _h}\left( \frac x{\varepsilon _h}\right) \right\|
_{L^p\left( \Omega ,\mathbf{R}^n\right) }\!\!\!.
\ea\hspace{-15mm}
\ee
As in the proof of Theorem~5 we have that
\[
\lim_{k\rightarrow \infty }\lim_{h\rightarrow \infty }\left\|
Du_h-Du_h^k\right\| _{L^p\left( \Omega ,\mathbf{R}^n\right) }=0,
\]
and, by Theorem~7, this implies that
\[
\lim_{k\rightarrow \infty }\lim_{h\rightarrow \infty }\left\|
Du_h^k-M_hDu_{*}^k-\sum_{i\in I_k}\chi _{\Omega _i^k}\left( x\right)
Dv^{M_hDu_{*}^k,x_i^k}\left( \frac x{\varepsilon _h}\right) \right\|
_{L^p\left( \Omega ,\mathbf{R}^n\right) }=0.
\]
This means that the theorem would be proved if we prove that
$\lim\limits_{k\rightarrow \infty }\lim\limits_{h\rightarrow \infty }$ acting on the last
term in (\ref{equ1thconcor}) is equal to~$0$. In order to prove this fact we
first make the following elementary estimates:
\be \label{equ2thconcor}
\ba{l}
\ds\left\| M_hDu_{*}^k+\sum\limits_{i\in I_k}\chi _{\Omega _i^k}\left(
x\right) Dv^{M_hDu_{*}^k,x_i^k}\left( \frac x{\varepsilon _h}\right)
-M_hDu-Dv^{M_hDu,\gamma _h}\left( \frac x{\varepsilon _h}\right) \right\|
_{L^p\left( \Omega ,\mathbf{R}^n\right) }^p
\vspace{3mm}\\
\ds \qquad =\sum\limits_{i\in I_k}\int_{\Omega _i^k}\left|
M_hDu_{*}^k+Dv^{M_hDu_{*}^k,x_i^k}\left( \frac x{\varepsilon _h}\right)
-M_hDu-Dv^{M_hDu,\gamma _h}\left( \frac x{\varepsilon _h}\right) \right|^p dx
\vspace{3mm}\\
\ds \qquad \leq C\sum\limits_{i\in I_k}\int_{\Omega _i^k}\left|
M_hDu_{*}^k+Dv^{M_hDu_{*}^k,x_i^k}\left( \frac x{\varepsilon _h}\right)
-M_hDu-Dv^{M_hDu,x_i^k}\left( \frac x{\varepsilon _h}\right) \right| ^pdx
\vspace{3mm}\\
\ds \qquad +C\sum\limits_{i\in I_k}\int_{\Omega _i^k}\left| Dv^{M_hDu,x_i^k}\left(
\frac x{\varepsilon _h}\right) -Dv^{M_hDu,\gamma _h}\left( \frac
x{\varepsilon _h}\right) \right| ^p\,dx.
\ea\hspace{-15mm}
\ee
We will study the two terms on the right hand side of (\ref{equ2thconcor})
separately, but first we define
\[
\xi _{h,*}^{j,k}=\frac 1{\left| Y_h^j\right| }\int_{Y_h^j}Du_{*}^k\,dx.
\]
By using a change of variables in (\ref{dalmasoineq}) and H\"{o}lder's
inequality, we find that we for the first term in (\ref{equ2thconcor}) have
the following estimate
\be\label{equ3thconcor}
\ba{l}
\ds \int_{\Omega _i^k}\left| M_hDu_{*}^k+Dv^{M_hDu_{*}^k,x_i^k}\left( \frac
x{\varepsilon _h}\right) -M_hDu-Dv^{M_hDu,x_i^k}\left( \frac x{\varepsilon
_h}\right) \right| ^pdx
\vspace{3mm}\\
\ds \qquad
\leq \sum\limits_{j\in J_h^{i,k}}\int_{Y_h^j}\left| \xi
_{h,*}^{j,k}+Dv^{\xi _{h,*}^{j,k},x_i^k}\left( \frac x{\varepsilon
_h}\right) -\xi _h^j-Dv^{\xi _h^j,x_i^k}\left( \frac x{\varepsilon
_h}\right) \right| ^p dx
\vspace{3mm}\\
\ds \qquad +\sum\limits_{j\in B_h^{i,k}}\int_{Y_h^j}\left| \xi _{h,*}^{j,k}+Dv^{\xi
_{h,*}^{j,k},x_i^k}\left( \frac x{\varepsilon _h}\right) -\xi _h^j-Dv^{\xi
_h^j,x_i^k}\left( \frac x{\varepsilon _h}\right) \right| ^p dx
\vspace{3mm}\\
\ds \qquad
\leq \sum\limits_{j\in J_h^{i,k}}C\left( 1+\left| \xi _{h,*}^{j,k}\right|
^p+\left| \xi _h^j\right| ^p\right) ^{\frac{\beta -\alpha -1}{\beta -\alpha
}}\left| \xi _{h,*}^{j,k}-\xi _h^j\right| ^{\frac p{\beta -\alpha }}\left|
Y_h^j\right|
\vspace{3mm}\\
\ds \qquad  +\sum\limits_{j\in B_h^{i,k}}C\left( 1+\left| \xi _{h,*}^{j,k}\right|
^p+\left| \xi _h^j\right| ^p\right) ^{\frac{\beta -\alpha -1}{\beta -\alpha
}}\left| \xi _{h,*}^{j,k}-\xi _h^j\right| ^{\frac p{\beta -\alpha }}\left|
Y_h^j\right|
\vspace{3mm}\\
\ds \qquad =C\int_{\Omega _i^{k,h}}\left( 1+\left| M_hDu_{*}^k\right| ^p+\left|
M_hDu\right| ^p\right) ^{\frac{\beta -\alpha -1}{\beta -\alpha }}\left|
M_hDu_{*}^k-M_hDu\right| ^{\frac p{\beta -\alpha }} dx
\ea
\ee
\[
\ba{l}
\ds  +C\int_{F_i^{k,h}}\left( 1+\left| M_hDu_{*}^k\right| ^p+\left|
M_hDu\right| ^p\right) ^{\frac{\beta -\alpha -1}{\beta -\alpha }}\left|
M_hDu_{*}^k-M_hDu\right| ^{\frac p{\beta -\alpha }}dx
\vspace{3mm}\\
\ds \leq C\left( \int_{\Omega _i^{k,h}}\!\!\left( 1+\left| M_hDu_{*}^k\right|
+\left| M_hDu\right| \right) ^pdx\right) ^{\frac{\beta -\alpha -1}{\beta
-\alpha }}\!\! \left( \int_{\Omega _i^{k,h}}\left| M_hDu_{*}^k-M_hDu\right|
^p dx\right) ^{\frac 1{\beta -\alpha }}
\vspace{3mm}\\
\ds  +C\left( \int_{F_i^{k,h}}\!\!\left( 1+\left| M_hDu_{*}^k\right| +\left|
M_hDu\right| \right) ^pdx\right) ^{\frac{\beta -\alpha -1}{\beta -\alpha }%
}\!\!\left( \int_{F_i^{k,h}}\left| M_hDu_{*}^k-M_hDu\right| ^p
dx\right) ^{\frac 1{\beta -\alpha }}
\vspace{3mm}\\
\ds
\leq C\left( \int_{\Omega _i^k}\left( 1+\left| Du_{*}^k\right| +\left|
Du\right| \right) ^pdx\right) ^{\frac{\beta -\alpha -1}{\beta -\alpha }
}\!\!\left( \int_{\Omega _i^k}\left| Du_{*}^k-Du\right| ^p dx\right) ^{\frac
1{\beta -\alpha }}
\vspace{3mm}\\
\ds
+C\left( \int_{F_i^{k,h}}\left( 1+\left| Du_{*}^k\right| +\left| Du\right|
\right) ^pdx\right) ^{\frac{\beta -\alpha -1}{\beta -\alpha }}
\!\! \left(
\int_{F_i^{k,h}}\left| Du_{*}^k-Du\right| ^p\,dx\right) ^{\frac 1{\beta
-\alpha }},
\ea\!
\]
where we used Jensen's inequality in the last step. Moreover, by using
(\ref{eq53}) and Jensen's inequality, we obtain that
\be\label{equ4thconcor}
\ba{l}
\ds \int_{\Omega _i^k}\left| Dv^{M_hDu,x_i^k}\left( \frac x{\varepsilon
_h}\right) -Dv^{M_hDu,\gamma _h}\left( \frac x{\varepsilon _h}\right)
\right| ^pdx
\vspace{3mm}\\
\ds \qquad \leq \sum\limits_{j\in J_h^{i,k}}\int_{Y_h^j}\left| Dv^{\xi
_h^j,x_i^k}\left( \frac x{\varepsilon _h}\right) -Dv^{\xi _h^j,x_h^j}\left(
\frac x{\varepsilon _h}\right) \right| ^p dx
\vspace{3mm}\\
\ds \qquad
+\sum\limits_{j\in B_h^{i,k}}\int_{Y_h^j}\left| Dv^{\xi _h^j,x_i^k}\left(
\frac x{\varepsilon _h}\right) -Dv^{\xi _h^j,x_h^j}\left( \frac
x{\varepsilon _h}\right) \right| ^pdx
\vspace{3mm}\\
\ds \qquad
\leq \sum\limits_{j\in J_h^{i,k}}C\omega \left(\frac 1k\right)\left( 1+\left| \xi
_h^j\right| \right) ^p\left| Y_h^j\right| +\sum\limits_{j\in
B_h^{i,k}}C\omega \left(\frac 1k+n\varepsilon _h\right)\left( 1+\left| \xi _h^j\right|
\right) ^p\left| Y_h^j\right|
\vspace{3mm}\\
\ds \qquad
=\omega \left(\frac 1k\right)C\int_{\Omega _i^{k,h}}\left( 1+\left| M_hDu\right|
^p\right) dx+\omega \left(\frac 1k+n\varepsilon _h\right)C\int_{F_i^{k,h}}\left(
1+\left| M_hDu\right| ^p\right) dx
\vspace{3mm}\\
\ds \qquad
\leq \omega \left(\frac 1k\right)C\int_{\Omega _i^k}\,\left( 1+\left| Du\right|
^p\right) dx+\omega \left(\frac 1k+n\varepsilon _h\right)C\int_{F_i^{k,h}}\left(
1+\left| Du\right| ^p\right) dx.
\ea\hspace{-10mm}
\ee
By combining (\ref{equ2thconcor}), (\ref{equ3thconcor}) and
(\ref{equ4thconcor}) we obtain that
\[
\ba{l}
\ds \left\| M_hDu_{*}^k+\sum\limits_{i\in I_k}\chi _{\Omega
_i^k}Dv^{M_hDu_{*}^k,x_i^k}\left( \frac x{\varepsilon _h}\right)
-M_hDu-Dv^{M_hDu,\gamma _h}\left( \frac x{\varepsilon _h}\right) \right\|
_{L^p\left( \Omega ,\mathbf{R}^n\right) }^p
\vspace{3mm}\\
\ds \qquad \leq C\left( \int_\Omega \left( 1+\left| Du_{*}^k\right| +\left| Du\right|
\right) ^pdx\right) ^{\frac{\beta -\alpha -1}{\beta -a}}\left( \int_\Omega
\left| Du_{*}^k-Du\right| ^pdx\right) ^{\frac 1{\beta -\alpha }}
\ea
\]
\[
\ba{l} \ds \qquad +C\sum\limits_{i\in I_k}\left(
\int_{F_i^{k,h}}\left( 1+\left| Du_{*}^k\right| +\left| Du\right|
\right) ^pdx\right) ^{\frac{\beta -\alpha -1}{\beta -a}}\left(
\int_{F_i^{k,h}}\left| Du_{*}^k-Du\right| ^p dx\right) ^{\frac
1{\beta -\alpha }}
\vspace{3mm}\\
\ds \qquad +C\omega \left(\frac 1k+n\varepsilon _h
\right)\left( \int_\Omega \left( 1+\left|
Du\right| ^p\right) dx+\sum\limits_{i\in I_k}\int_{F_i^{k,h}} \left(
1+\left| Du\right| ^p\right) dx\right) .
\ea
\]
Moreover, by noting that $\left| F_i^{k,h}\right| \rightarrow 0$ as
$h\rightarrow \infty ,$ $\left\| Du\right\| _{L^p\left( \Omega ,
\mathbf{R}^n\right) }\leq C,$ $\left\| Du_{*}^k\right\| _{L^p\left( \Omega ,
\mathbf{R}^n\right) }\leq C,$ and taking (\ref{eq44}) into account, we obtain that
\[
\ba{l}
\ds \lim_{k\rightarrow \infty }\lim_{h\rightarrow \infty }\Biggl\|
M_hDu_{*}^k+\sum_{i\in I_k}\chi _{\Omega _i^k}Dv^{M_hDu_{*}^k,x_i^k}\left(
\frac x{\varepsilon _h}\right)
\vspace{3mm}\\
\ds \qquad \qquad  -M_hDu-Dv^{M_hDu,\gamma _h}\left( \frac
x{\varepsilon _h}\right) \Biggr\| _{L^p\left( \Omega ,\mathbf{R}^n\right) }=0,
\ea
\]
and the theorem is proved.\hfill \rule{3mm}{3mm}

\subsection*{Acknowledgments}

The author wishes to thank Professor
Lars-Erik Persson and Dr. Peter Wall for some~help\-ful comments and
suggestions which have improved this paper.

\label{bystrom-lastpage}

\end{document}